\documentstyle[fleqn,12pt]{article}
\begin{document}
\pagenumbering{arabic}\setcounter{page}{1}\pagestyle{plain}
\baselineskip=16pt

\thispagestyle{empty}
\rightline{MSUMB 96-02, October 1996} 
\vspace{1.4cm}

\begin{center}
{\bf $h$-deformation of grassmann supergroup $Gr(1|1)$ } \end{center}

\vspace{1cm}
Salih \c Celik \footnote {E-mail: celik@msu.edu.tr}
\footnote {New E-mail: sacelik@yildiz.edu.tr}\\ 
Mimar Sinan University, Department of Mathematics, \\
80690 Besiktas, Istanbul, TURKEY.

\vspace{2.5cm}
{\bf Abstract}

We introduce the $h$-deformation of the algebra of functions on the 
grassmann supergroup Gr$(1|1)$ via a contraction of Gr$_q(1|1)$. 

\vfill\eject
Recently matrix groups like GL$(2)$, GL$(1|1)$, {\it etc.}, 
were generalized in two ways. Both are based on the deformation 
of the algebra of functions on the groups generated by 
coordinate functions $t^i{}_j$ that commute 
$$[t^i{}_j, t^k{}_l] = 0.$$

In the $q$-deformation of matrix groups [1], these 
commutation relations are determined by a matrix $R_q$ 
so that the functions do not commute but satisfy the equation 
$$R_q (T \otimes T) = (I \otimes T)(T \otimes I) R_q.$$
In this equation the elements of $R_q$ are numbers but 
the matrix $T = (t^i{}_j)$ is formed by generally non-commuting 
elements of an algebra. 

Another type of deformation, so called the "$h$-deformation", which 
is a new class of quantum deformations of Lie groups and Lie algebras 
has recently been intensively studied [2-4]. These deformations 
can be important since the deformation parameter $h$ can directly be 
identified with the Planck constant, which makes them interesting 
for physical applications. The other important point is that the 
standard $q$-deformation and the $h$-deformation have been characterized 
as the only two deformations of SL(2) [3]. 

The goal of this paper is to develop the $h$-deformation for the case 
of grassmann matrix supergroups. We study the simplest grassmann 
supergroup Gr$(1|1)$, here. The one parametric $q$-deformation of Gr$(1|1)$, 
Gr$_q(1|1)$, was given in Ref. 5. Before discussing the $h$-deformation 
of Gr$(1|1)$, we give some notations and useful formulas in following. 

A grassmann supermatrix $\hat{T}$, that is, an element of 
Gr$(1|1)$, is of the form 
$$ \hat{T} = \left(\matrix{ \alpha & b \cr c & \delta \cr} \right)$$
with two odd (greek letters) and two even (latin letters) matrix 
elements. The following closely follows the approach of Ref. 5. 

In this paper we denote $q$-deformed objects by primed quantities. Unprimed 
quantities represent transformed coordinates. 

Let us consider the following quantum superplane and its dual, 
denoted by $A_q$ and $A_q^*$, respectively [6]: 
$$U' = \left (\begin{array}{c} x' \\ \xi' \end{array} \right ) \in A_q ~
      \Longleftrightarrow~ x' \xi' - q \xi' x' = 0,~ \xi'^2 = 0, \eqno(1) $$
and its dual 
$$\hat{U}' = \left (\begin{array}{c} \eta' \\ y' \end{array} \right ) 
   \in A_q^* ~\Longleftrightarrow~ 
   \eta'^2 = 0,~ \eta' y' - q^{-1} y' \eta' = 0.    \eqno(2) $$

Suppose that the matrix elements of $\hat{T}'$ (anti-)commute 
with the coordinates of $A_q$ and $A_q^*$. Then, the endomorphisms 
$$ \hat{T}' : A_q \longrightarrow A_q^* ~~\mbox{and}~~ 
   \hat{T}' : A_q^* \longrightarrow A_q \eqno(3)$$
impose the following bilinear product relations among the matrix 
elements of $\hat{T}'$:
$$ \alpha' b' = q^{-1} b' \alpha', \quad \alpha' c' = q^{-1} c' \alpha', $$
$$ \delta' b' = q^{-1} b' \delta', \quad \delta' c' = q^{-1} c' \delta', $$
$$\alpha' \delta' + \delta' \alpha' = 0, \quad \alpha'^2 = 0 = \delta'^2, 
   \eqno(4)$$
$$ b'c' = c'b' + (q - q^{-1}) \delta' \alpha' $$
where $q$ is a non-zero complex number and $q^2 \neq 1$. 
These relations are the defining relations of Gr$_q(1|1)$. 

We introduce new coordinates $x$ and $\xi$ by 
$$U = g^{-1} U' \eqno(5) $$
where 
$$ g = \left(\matrix{1 & 0 \cr h/(q - 1) & 1 \cr} \right) \eqno(6) $$
as in Ref. 7. Here the deformation parameter $h$ is a grassmann 
number which has the following properties 
$$h^2 = 0 ~~\mbox{and}~~ h \xi = - \xi h. $$
Now, in the limit $q \rightarrow 1$ we get the following exchange 
relations, which define the $h$-superplane $A_h$ 
$$x \xi = \xi x + h x^2, \quad \xi^2 = - h x \xi. \eqno(7)$$

Similarly, using the relations (2) with (6), one has [after taking the 
limit $q \rightarrow 1$] the dual $h$-superplane $A^*_h$ as generated 
by $\eta$ and $y$ with the relations 
$$\eta^2 = 0, \quad \eta y = y \eta. \eqno(8)$$

We now consider the endomorphisms 
$$ \hat{T} : A_h \longrightarrow A_h^* ~~\mbox{and}~~ 
   \hat{T} : A_h^* \longrightarrow A_h. \eqno(9)$$
Then, we define the corresponding $h$-deformation of the grassmann 
supergroup Gr$(1|1)$ as a quantum matrix supergroup Gr$_h(1|1)$ generated 
by $\alpha$, $b$, $c$ and $\delta$ which satisfy the following 
$h$-commutation relations 
$$\alpha b = b \alpha + h b^2, \quad 
   \alpha c = c \alpha + h(bc + \delta \alpha), $$
$$ \delta b = b \delta - h b^2, \quad 
   \delta c = c \delta - h(cb - \alpha \delta), $$
$$ \alpha^2 = h \alpha b, \quad 
   \alpha \delta = - \delta \alpha + h(\delta b - b \alpha), \eqno(10)$$
$$ \delta^2 = - h \delta b, \quad 
   bc = cb + h (b \delta + \alpha b). $$

It can be check that the maps 
$$\mu : A_h \longrightarrow Gr(1|1) \otimes A^*_h, \quad 
  \mu^* : A^*_h \longrightarrow Gr(1|1) \otimes A_h \eqno(11)$$
such that 
$$ \mu(U) = \hat{T} \otimes \hat{U}, ~~i.e. ~~ 
   \mu(u_i) = t^j_i \otimes \hat{u}_j, $$
$$ \mu^*(\hat{U}) = \hat{T} \otimes U, ~~i.e. ~~ 
   \mu^*(\hat{u}_i) = t^j_i \otimes u_j \eqno(12)$$
define the co-action of the quantum group $Gr_h(1|1)$ on the quantum 
superplane $A_h$ and its dual $A_h^*$, respectively. 

An interesting feature is that the second column of the matrix 
$\hat{T} \in$ Gr$_h(1|1)$ is isomorphic to $A_h$, though this is 
not so in GL$_h(1|1)$ [7]. 

Alternativelly, the relations (10) can be obtained using the following 
similarity transformation which given by Aghamohammadi {\it et al} [4]: 
$$\hat{T}' = g \hat{T} g^{-1}. \eqno(13)$$
To do this, we use the relations (4) and next the limit $q \rightarrow 1$. 

The quantum (dual) superdeterminant of $\hat{T}$ is again defined as [5] 
$$ \hat{D}_h = bc^{-1} - \alpha c^{-1} \delta c^{-1} 
                 = c^{-1}b - c^{-1} \alpha c^{-1} \delta \eqno(14)$$
which is independent of the relations (10). It is easy to verify that 
$\hat{D}_h$ commutes with all matrix elements of $\hat{T}$ ~(and $h$), 
that is, $\hat{D}_h$ belongs to the centre of the algebra 
$$\hat{T} \hat{D} = \hat{D} \hat{T}. $$

The inverse of $\hat{T}$ is defined as [5] 
$$\hat{T}^{-1}  
 = \left(\matrix{- c^{-1} \delta b^{-1} & c^{-1} + 
                             c^{-1} \delta b^{-1} \alpha c^{-1} \cr \cr
                  b^{-1} + b^{-1} \alpha c^{-1} \delta b^{-1} & 
                               - b^{-1} \alpha c^{-1} \cr}
\right) \eqno(15) $$
provided $b$ and $c$ are invertible. 

We shall now obtain an R-matrix for the quantum grassmann 
supergroup Gr$_h(1|1)$ from the R-matrix of Gr$_q(1|1)$ which is 
introduced in Ref. 8. 

The associative algebra (4) is equivalent to 
$$R_q \hat{T}_1' \hat{T}_2' = - 
  \hat{T}_2' \hat{T}_1' R_q \eqno(16)$$
where [8] 
$$R_q 
 = \left(\matrix{ q + q^{-1} & 0          & 0          & 0 \cr 
                        0    & - 2        & q^{-1} - q & 0 \cr
                        0    & q - q^{-1} & - 2  & 0 \cr 
                        0    & 0          & 0    & q + q^{-1}  \cr }
\right). \eqno(17)$$
Here, we used the tensoring convention 
$$(\hat{T}_1)^{ij}{}_{kl} 
   = (\hat{T} \otimes I)^{ij}{}_{kl} 
   = (-1)^{k(j+l)} \hat{T}^i{}_k \delta^j{}_l, $$
$$ (\hat{T}_2)^{ij}{}_{kl} 
   = (I \otimes \hat{T})^{ij}{}_{kl} 
   = (-1)^{i(j+l)} \hat{T}^j{}_l \delta^i{}_k. \eqno(18)$$
Note that although the algebra (4) is an associative algebra of the 
matrix elements of $\hat{T}'$, $R_q$ does not satisfy the (graded) quantum 
Yang-Baxter equation (QYBE) 
$$R_{12}R_{13}R_{23} = R_{23}R_{13}R_{12}. \eqno(19) $$
Thus the QYBE is not a necessary condition for associativity. 
However, as interesting point, if we decompose the matrix $R_q$ is of the form 
$$R_q = R_q^1 + R_q^2 \eqno(20) $$
where 
$$R_q^1 
 = \left(\matrix{ q &   0         & 0    & 0 \cr 
                  0 & - 1         & 0    & 0 \cr
                  0 &  q - q^{-1} & - 1  & 0 \cr 
                  0 &   0         & 0    & q^{-1}  \cr }\right), \quad 
R_q^2 
 = \left(\matrix{ q^{-1} & 0   &  0         & 0 \cr 
                    0    & - 1 & q^{-1} - q & 0 \cr
                    0    & 0   & - 1        & 0 \cr 
                    0    & 0   &  0         & q \cr }
\right) \eqno(21)$$
then both matrices $R_q^1$ and $R_q^2$ do satisfy  the (graded) QYBE. 
Thus we can say that the addition of two matrices satisfying the (graded) 
QYBE does not need to satisfy the (graded) QYBE. 
Also the equation (16) can be written of the form 
$$R_q^1 T_1' T_2' = - T_2' T_1' R_q^2. \eqno(22)$$

Let $g_i$ and $\hat{T}_i$ $(i = 1, 2)$ be graded in the 
same manner as $\hat{T}_i'$. 
It is clear that 
$$\hat{T}_1' \hat{T}_2' = 
  G_h \hat{T}_1 \hat{T}_2 G_h ~~\mbox{and}~~ 
  \hat{T}_2'\hat{T}_1' = 
  G_{-h} \hat{T}_2 \hat{T}_1 G_{-h}, \eqno(23)$$
where 
$$G_h = g_1 g_2^{-1}, \quad G_{-h} = G_h^{-1}. \eqno(24)$$
Substituting (23) into (16) we arrive at the equation 
$$G_h R_q G_h \hat{T}_1 \hat{T}_2 = - 
  \hat{T}_2 \hat{T}_1 G_{-h} R_q G_{-h}. \eqno(25)$$
If we define the R-matrix $R_h$ as 
$$R_h = \lim_{q \rightarrow 1} \big(G_hR_qG_h\big), \eqno(26)$$
[after dividing by 2] we get the following R-matrix $R_h$ 
$$R_h 
 = \left(\matrix{ 1 &   0  &   0  & 0 \cr 
                 -h & - 1  &   0  & 0 \cr
                 -h &   0  & - 1  & 0 \cr 
                  0 &   h  & - h  & 1   \cr }
\right) \eqno(27)$$
which gives the relations (10) with the equation 
$$R_h \hat{T}_1 \hat{T}_2 = - 
  \hat{T}_2 \hat{T}_1 R_{-h}. \eqno(28)$$
Note that although $R_h$ does not satisfy the usual graded and ungraded 
YBE (19) it satisfies the equation 
$$(R_{12}R_{13}R_{23})(h) = (R_{23}R_{13}R_{12})(- h), \eqno(29) $$
or, equivalently, 
$$R_{12}R_{13}R_{23} + R_{23}R_{13}R_{12} = 2I_8, \eqno(30) $$
where $I_8$ is the 8x8 unit matrix and
$$R_{12} = R \otimes I_2, ~~etc. $$
Since 
$$R_h^2 = I,$$
the R-matrix $R_h$ has two eigenvalues $\pm 1$. 

{\bf Acknowledgement}

This work was supported in part by {\bf T. B. T. A. K.} the 
Turkish Scientific and Technical Research Council. 

\def\refname{References}


\begin{thebibliography}{99}
\bibitem{D:gnus} V. G. Drinfeld, {\em Quantum groups}, in {\it Proc. } IMS, Berkeley, 1986; 

N. Y. Reshetikhin, L. A. Takhtajan and L. D. Faddeev, 
    {\it Leningrad Math. J. } {\bf 1}, 193 (1990). 
\bibitem{DMZ:gnus} 
E. E. Demidov, Yu. I. Manin, E. E. Mukhin and D. V. Zhdanovich, 
     {\it Prog. Theor. Phys. Suppl.} {\bf 102}, 203 (1990); 

H. Ewen, O. Ogievetsky and J. Wess, 
    {\it Lett. Math. Phys.} {\bf 22}, 297 (1991); 

S. Zakrzewski, {\it Lett. Math. Phys.} {\bf 22}, 287 (1991);

S. L. Woronowicz, {\it Rep. Math. Phys.} {\bf 30}, 259 (1991);

C. H. Ohn, {\it Lett. Math. Phys.} {\bf 25}, 85 (1992). 

\bibitem{K:gnus} 
B. A. Kupershmidt, {\it J. Phys. A: Math. Gen.} {\bf 25}, L1239 (1992); 

V. Karimipour, {\it Lett. Math. Phys.} {\bf 30}, 87 (1994); 

A. Aghamohammadi, {\it Mod. Phys. Lett. A} {\bf 8}, 2607 (1993).

S. Celik and E. Hizel, {\it $h$-deformation of $Gr(2)$}, to 
appear in {\it J. Phys. A: Math. Gen.} {\bf 30} (1997), 4677-4680. 
\bibitem{AKS:gnus} 
A. Aghamohammadi, M. Khorrami and A. Shariati, 
     {\it J. Phys. A: Math. Gen.} {\bf 28}, L225 (1995).
\bibitem{SS:gnus} 
S. Celik and S. A. Celik, {\it Balkan Phys. Lett.} {\bf 3}, 188 (1995). 
\bibitem{M:gnus} 
Yu. I. Manin, {\em Commun. Math. Phys.} {\bf 123}, 163 (1989). 
\bibitem{DP:gnus} 
L. Dabrowski and P. Parashar, {\em $h$-deformation of GL$(1|1)$} 
     to appear in {\it Lett.} {Math. Phys.} 
\bibitem{SC:gnus} 
S. Celik, {\it J. Math. Phys.} {\bf 37}, 3568 (1996). 

\end{thebibliography}
\end{document}